\documentclass[12pt]{article}
\usepackage{amssymb}
\usepackage{amsfonts}
\usepackage{times}
\usepackage{mathptmx}
\usepackage{amsmath}
\usepackage[usenames]{color}
\usepackage{mathrsfs}
\usepackage{amsfonts}
\usepackage{amssymb,amsmath}
\usepackage{CJK}
\usepackage{cite}
\usepackage{cases}
\usepackage{amsthm}
\def\cl{\centerline}

\def\vs{\vspace*}

\def\C{\mathbb{C}}

\numberwithin{equation}{section}
\newtheorem{theo}{Theorem}[section]
\newtheorem{defi}[theo]{Definition}
\newtheorem{coro}[theo]{Corollary}

\newtheorem{remark}[theo]{Remark}
\newtheorem{examp}[theo]{Example}
\newtheorem{prop}[theo]{Proposition}

\begin{document}
\begin{center}
{\bf\large Quadratic Leibniz conformal algebras \,{$^*\,$}}
\footnote {\,$^*\,$ Supported by the Zhejiang Provincial Natural Science Foundation of China (No.~LQ16A010011) and the National Natural Science Foundation of China (No.~11501515, 11871421).

$^{\,\dag}$Corresponding author: Y.~Hong.
}
\end{center}

\cl{Jinsen Zhou$^{1}$, Yanyong Hong$^{2,\dag}$
}

\cl{\small $^{1}$zjs9932@126.com}
\cl{\small $^{1}$School of Information Engineering, Longyan University, Longyan 364012, Fujian, P. R. China}
\cl{\small $^{2}$hongyanyong2008@yahoo.com}
\cl{\small $^{2}$College of Science, Zhejiang Agriculture and Forestry University,}
\cl{\small Hangzhou 311300, Zhejiang, P. R. China}

\vs{8pt}

{\small\footnotesize
\parskip .005 truein
\baselineskip 3pt \lineskip 3pt
\noindent{{\bf Abstract:} In this paper, we study a class of Leibniz conformal algebras called quadratic Leibniz conformal algebras. An equivalent characterization of a Leibniz conformal algebra $R=\mathbb{C}[\partial]V$ through three algebraic operations on $V$ are given. By this characterization, several constructions of quadratic Leibniz conformal algebras are presented. Moreover, one-dimensional central extensions of
quadratic Leibniz conformal algebras are considered using some bilinear forms on $V$. In particular, we also study one-dimensional Leibniz central extensions of quadratic Lie conformal algebras.
\vs{5pt}


}}
\parskip .001 truein\baselineskip 6pt \lineskip 6pt

\section{Introduction}
Throughout this paper, we denote by $\C$ the set of complex numbers,
$\mathbb{Z}$ the set of integer
numbers and $\mathbb{Z}_{+}$ the set of non-negative integer numbers. All vector spaces are over $\mathbb{C}$ and tensors over $\mathbb{C}$ are denoted by  $\otimes$.
Moreover, if $V$ is a vector space, then the space of polynomials of $\lambda$ with coefficients in $V$ is denoted by $V[\lambda]$.

Lie conformal algebra, introduced by Kac in \cite{K1}, gives an axiomatic description of the singular part of the operator product expansion of chiral fields in conformal field theory. It is an useful tool to study vertex algebras (see \cite{K1}) and has many applications in the theory of infinite-dimensional Lie algebras. Moreover, Lie conformal algebras have close connections to Hamiltonian formalism in the theory of nonlinear evolution equations (see \cite{BDK}). In fact, a Lie conformal algebra is a Lie pseudo-algebra which can be seen a Lie algebra in a pseudo-tensor category (see \cite{BDK1}). Structure theory and representation theory of finite Lie conformal algebras which are finitely generated as $\mathbb{C}[\partial]$-modules are well developed (see \cite{CK1, DK, CK2, BKV} and so on).

Conformal algebras
are quite intriguing subjects in the purely algebraic viewpoint. One can define the conformal analogue
of a variety of ``usual" algebras such as Lie conformal algebras,
associative conformal algebras, etc. The theory of
conformal algebras sheds new light on the problem of classification
of infinite-dimensional algebras of the corresponding ``classical"
variety.

Leibniz algebras are the non-commutative analogs of Lie algebras. Leibniz algebra, first introduced by Bloh in \cite{B}, and reintroduced by Loday in \cite{L}, arose naturally during their studying on periodicity phenomena in algebraic $K$-theory. The name of left (right) Leibniz algebra comes from that the left (right) multiplication is a derivation. We call the left Leibniz algebra as Leibniz algebra in this paper. Note that Leibniz algebras are di-Lie algebras. There is a replication procedure described in \cite{gk1} which allows to get the defining identities of a class of di-algebras based on a given variety of algebras.

The topic of this paper is about Leibniz conformal algebra. Leibniz conformal algebra was introduced in \cite{BKV} and its cohomology theory was investigated in \cite{BKV} and \cite{Z}. In addition, Leibniz pseudoalgebra was studied in \cite{Wu}. But, there are few examples of Leibniz conformal algebras which are not Lie conformal algebras. In this paper, our aim is to provide an efficient way to construct Leibniz conformal algebras which are not Lie conformal algebras.
By the correspondence of Leibniz conformal algebras and (infinite-dimensional) Leibniz algebras, it is also useful to construct infinite-dimensional Leibniz algebras which are not Lie algebras. Therefore, it is of signification and interesting. For Lie conformal algebras, as we know, most of known examples are quadratic Lie conformal algebras which were studied in \cite{GD, Xu1}. It was essentially stated in \cite{GD} that a ``quadratic" Lie conformal algebra is equivalent to a vector space with two algebra structures. One is a Lie algebra structure, the other is a Novikov algebra structure and they satisfy a compatibility condition. This algebra structure is called a Gel'fand-Dorfman bialgebra by Xu in \cite{Xu1}. For not confusing the concept of bialgebra in Hopf algebra, we call it Gel'fand-Dorfman algebra in short. In fact, a quadratic Lie conformal algebra corresponds to a Hamiltonian pair in \cite{GD}, which plays fundamental roles in completely integrable systems. Moreover, several constructions of Gel'fand-Dorfman algebras were presented by Xu in \cite{Xu1}. Therefore, it is an useful method to construct
Lie conformal algebras by using Gel'fand-Dorfman algebras. Note that Lie conformal algebras are Leibniz
conformal algebras. Motivated by this, in this paper, we investigate what are the algebra structures underlying quadratic Leibniz conformal algebras. Since Gel'fand-Dorfman algebras correspond to
Lie conformal algebras, if we find other algebra structures which are not Gel'fand-Dorfman algebras underlying Leibniz conformal algebras, we can construct Leibniz conformal algebras which are not Lie conformal algebras.
We show that except Gel'fand-Dorfman algebra, there is also an algebra structure called Perm-Leibniz algebra underlying quadratic Leibniz conformal algebra. Through Perm-Leibniz algebras, we can construct many Leibniz conformal algebras which are not Lie conformal algebras. Moreover, one-dimensional central extensions of
quadratic Leibniz conformal algebra $R=\mathbb{C}[\partial]V$ are considered using some bilinear forms on $V$. In particular, we also study one-dimensional Leibniz central extensions of quadratic Lie conformal algebras. It should be pointed out that when we consider the case of superalgebra,  similar results can be directly obtained from the ones in this paper.

This paper is organized as follows. In Section 2, some basic definitions and some facts about Lie conformal algebras and Leibniz conformal algebras are recalled. We also introduce the definitions of quadratic Lie and Leibniz conformal algebras. In Section 3, an equivalent characterization of quadratic Leibniz conformal algebra is given. Through this equivalent characterization, we find that a vector space with two algebra structures called Perm-Leibniz algebra can be used to construct Leibniz conformal algebras which are not Lie conformal algebras. Moreover, several constructions and examples are given. In Section 4, we investigate  one-dimensional central extensions of
quadratic Leibniz conformal algebras $R=\mathbb{C}[\partial]V$ corresponding to Perm-Leibniz algebras using some bilinear forms on $V$. In particular, one-dimensional Leibniz central extensions of quadratic Lie conformal algebras are also considered.

\section{Preliminaries}
In this section, we will introduce some basic definitions and some facts
about Lie conformal algebras and Leibniz conformal algebras.

First, let us recall the definition of Leibniz algebra.
\begin{defi}
A \emph{Leibniz algebra} $L$ is a vector space $L$ with a bilinear operation $[\cdot, \cdot]$ satisfying the following axioms:
\begin{eqnarray}\label{eq1}
[a,[b,c]]=[[a,b],c]+[b,[a,c]],
\end{eqnarray}
for $a$, $b$, $c\in L$.
\end{defi}
\begin{remark}
(\ref{eq1}) is called \emph{left Leibniz identity}. Similarly, there is
a \emph{right Leibniz identity}:
\begin{eqnarray}
[a,[b,c]]=[[a,b],c]-[[a,c],b],
\end{eqnarray}
for $a$, $b$, $c\in L$. If a vector space $L$ with a bilinear operation $[\cdot, \cdot]$ satisfying the right Leibniz identity, then $L$ is called a \emph{right Leibniz algebra}.

Note that, any Leibniz algebra $(L,[\cdot,\cdot])$ can be endowed with a right Leibniz algebra structure $(L,[\cdot,\cdot]^{'})$ via $[a,b]^{'}=-[b,a]$
for any $a$, $b\in L$ and vice versa.
\end{remark}

Obviously, all Lie algebras are Leibniz algebras. Next, we consider
the conformal versions of these algebras.

\begin{defi}(see \cite{K1}) \label{250}\rm
A \emph{conformal algebra} $R$ is a $\mathbb{C}[\partial]$-module $R$ endowed with a $\C$-bilinear map
\begin{equation*}\label{251}
R\times R\rightarrow R[\lambda], ~~~~~~a\times b\mapsto [a{}\, _\lambda \, b],
\end{equation*}
satisfying the following axiom~($a, b\in R$):
\begin{equation*}\label{251}
\aligned
Conformal~~sesquilinearity:~~~~[\partial a\,{}_\lambda \,b]=-\lambda[a\,{}_\lambda\, b],\ \ \ \
[a\,{}_\lambda \,\partial b]=(\partial+\lambda)[a\,{}_\lambda\, b].
\endaligned
\end{equation*}

A \emph{Lie conformal algebra} $(R,[\cdot_{\lambda}\cdot])$ is a conformal algebra and satisfies the following axioms ~($a, b, c\in R$):
\begin{equation*}
\aligned
&Skew~~symmetry:~~~~[a\, {}_\lambda\, b]=-[b\,{}_{-\lambda-\partial}\,a];\\
&Jacobi~~identity:~~~~[a\,{}_\lambda\,[b\,{}_\mu\, c]]=[[a\,{}_\lambda\, b]\,{}_{\lambda+\mu}\, c]+[b\,{}_\mu\,[a\,{}_\lambda \,c]].
\endaligned
\end{equation*}
\end{defi}

\begin{defi}(see \cite{BKV}) \label{252}\rm
A \emph{Leibniz conformal algebra} $(R,[\cdot_\lambda \cdot])$ is a conformal algebra and satisfies
\begin{equation*}
Leibniz~~identity:~~~~[a\,{}_\lambda\,[b\,{}_\mu\, c]]=[[a\,{}_\lambda\, b]\,{}_{\lambda+\mu}\, c]+[b\,{}_\mu\,[a\,{}_\lambda \,c]],
\end{equation*}
for any $a$, $b$, $c\in R$.
\end{defi}

A Lie or Leibniz conformal algebra $R$ is called \emph{finite}, if it is finitely generated as a
$\mathbb{C}[\partial]$-module; otherwise, it is said to be \emph{infinite}.

\begin{remark}
Similar to the classical case, there is the definition of right Leibniz conformal algebra. A \emph{right Leibniz conformal algebra} $(R,[\cdot_\lambda \cdot])$ is a conformal algebra and satisfies the right Leibniz identity
\begin{eqnarray*}
[a\,{}_\lambda\,[b\,{}_\mu\, c]]=[[a\,{}_\lambda\, b]\,{}_{\lambda+\mu}\, c]-[[a\,{}_\lambda \,c]\,{}_{-\mu-\partial}\, b],
\end{eqnarray*}
where $a$, $b$ and $c\in R$.

Similarly, it is easy to check that any Leibniz conformal algebra $(R,[\cdot_\lambda \cdot])$ can be endowed with a right Leibniz conformal algebra structure $(R, [\cdot_\lambda \cdot]^{'})$ via
$[a_\lambda b]^{'}=-[b_{-\lambda-\partial} a]$ for
any $a$ and $b\in R$ and vice versa. Since the two algebra structures can be obtained from each other, in this paper, we only study Leibniz conformal algebras.
\end{remark}

\begin{remark}
Obviously, all Lie conformal algebras are Leibniz conformal algebras.
\end{remark}

\begin{examp}
Let $L$ be a Leibniz algebra. Then $\text{Cur}L=\mathbb{C}[\partial]\otimes L$ can be endowed with a Lebiniz conformal algebra structure as follows
\begin{eqnarray}
[a_\lambda b]=[a,b],~~~~~~a,~~b\in L.
\end{eqnarray}
$\text{Cur}L$ is called the \emph{current Leibniz conformal algebra} associated with
$L$.
\end{examp}

Suppose $R$ is a Lie (or Leibniz) conformal algebra. Write $[a_\lambda b]=\sum_{i\in \mathbb{Z}_+}\frac{\lambda^i}{i!} a_{(i)}b$, where $a_{(i)}b\in R$ for every $i$. There is a natural Lie (or Leibniz) algebra associated with it.
Let \text{Coeff}$(R)$ be the quotient
of the vector space with the basis $a_{n}$ $(a\in R, n\in\mathbb{Z})$ by
the subspace spanned over $\mathbb{C}$ by
elements:
$$(\alpha a)_{n}-\alpha a_{n},~~(a+b)_{n}-a_{n}-b_{n},~~(\partial
a)_{n}+na_{n-1},~~~\text{where}~~a,~~b\in R,~~\alpha\in \mathbb{C},~~n\in
\mathbb{Z}.$$
The operation on \text{Coeff}$(R)$ is defined as follows:
\begin{eqnarray}
[a_{m}, b_{n}]=\sum_{j\in \mathbb{Z}_{+}}\left(\begin{array}{ccc}
m\\j\end{array}\right)(a_{(j)}b)_{m+n-j}.\end{eqnarray}
Then \text{Coeff}$(R)$ is a Lie (or Leibniz) algebra. Of course,
if $R$ is not torsion as a $\mathbb{C}[\partial]$-module, then \text{Coeff}$(R)$ is infinite-dimensional.

\begin{defi}\label{254}\rm
Let $R$ be a conformal algebra. If there exists some vector space $V$ such that $R=\C[\partial]V$ is a free $\C[\partial]$-module over $V$ and for all $a$, $b\in V$, the corresponding $\lambda$-product is of the following form:
\begin{eqnarray}
[a_\lambda b]=\partial u+\lambda v+w, ~~~~~\text{$u$, $v$, $w\in V$,}
\end{eqnarray}
then $R$ is called \emph{quadratic}.
\end{defi}

Finally, let us recall the definitions of Novikov algebra and Gel'fand-Dorfman algebra.

\begin{defi}
A \emph{Novikov algebra} $(A, \circ)$ is a vector space $A$ with a bilinear operation
``$\circ$" satisfying the following axioms: for all $a$, $b$, $c\in A$,
\begin{eqnarray*}
&a\circ (b\circ c)=b\circ(a\circ c),\\
&(a\circ b)\circ c-a\circ (b\circ c)=(a\circ c)\circ b-a\circ(c\circ b).
\end{eqnarray*}
\end{defi}

\begin{remark}
Novikov algebra corresponds to a
certain Hamiltonian operator (see \cite{GD}). Such an algebraic structure is also related with Poisson structure of
hydrodynamic type (see \cite{BN}). The name ``Novikov algebra" was given by Osborn in
\cite{Os}.
\end{remark}

\begin{defi}
A \emph{Gel'fand-Dorfman algebra} $(A, \circ, [\cdot, \cdot] )$ is a vector space $A$ with two algebraic operations $[\cdot,\cdot]$ and $\circ$ such that $(A,[\cdot,\cdot])$ forms a Lie algebra, $(A,\circ)$ forms a Novikov algebra and the following compatibility condition holds
\begin{eqnarray}\label{xx1}
a\circ [b,c]+[a,b\circ c]=b\circ[a,c]+[b,a\circ c]+[a,b]\circ c,
\end{eqnarray}
for $a$, $b$, and $c\in A$.
\end{defi}
\begin{remark}
It should be pointed out that a similar algebra structure in \cite{Xu1} is called \emph{Gel'fand-Dorfman bialgebra}.  For not confusing the concept of bialgebra in Hopf algebra, we call it Gel'fand-Dorfman algebra in short.
\end{remark}
\begin{theo}\label{th1}(see \cite{Xu1})
A  quadratic conformal algebra $R=\mathbb{C}[\partial]V$ is Lie if and only if $(V,\circ, [\cdot,\cdot])$ is a Gel'fand-Dorfman algebra relative to
\begin{eqnarray}\label{xy1}
[a_\lambda b]=\partial(a\circ b)+\lambda(a\circ b+b\circ a)+[a,b],
\end{eqnarray}
for all $a$ and $b\in V$.
\end{theo}
\begin{remark}
In \cite{Xu1}, the correspondence between quadratic Lie conformal algebra and Gel'fand-Dorfman bialgebra is $[a_\lambda b]=\partial(b\circ a)+\lambda(a\circ b+b\circ a)+[b,a]$ which is a little different from (\ref{xy1}). Since (\ref{xy1}) makes all expressions much more natural, we choose it.
\end{remark}

\section{Characterization of quadratic Leibniz conformal algebras}
In this section, we give an equivalent characterization of quadratic Leibniz conformal algebra and present several constructions of quadratic Leibniz conformal algebras.

\begin{theo}\label{th2}
Let $V$ be a vector space. A quadratic Leibniz conformal algebra $R=\mathbb{C}[\partial]V$ is equivalent to the quadruple $(V, \vdash, \dashv, [\cdot,\cdot])$ where
$\vdash$, $\dashv$ are two bilinear operations on $V$, $[\cdot,\cdot]$ is a Leibniz
algebra operation, and they satisfy the following conditions:
\begin{gather}
\label{t1}x\vdash (y\vdash z)=y\vdash (x\vdash z),\\
\label{t2}x\vdash (y\dashv z)=y\dashv(x\dashv z),\\
\label{t3}x\vdash(y\vdash z-z\dashv y)=(x\vdash y)\vdash z-(x\vdash z)\dashv y,\\
\label{tt1}(x\vdash y-x\dashv y)\vdash z=0,\\
\label{t4}x\dashv (y\vdash z-y\dashv z)=0,\\
\label{t5}x\dashv(y\vdash z-z\dashv y)=(x\dashv y)\dashv z-(x\dashv z)\dashv y,\\
\label{t7}x\vdash[y,z]+[x,y\vdash z]=y\vdash[x,z]+[y,x\vdash z]+[x,y]\vdash z,\\
\label{t8}x\vdash [y,z]+[x,z]\dashv y=[x,z\dashv y]+[x\dashv y,z]+z\dashv [y,x],\\
\label{t9}[x\vdash y-x\dashv y,z]=z\dashv([x,y]+[y,x]),
\end{gather}
for any $x$, $y$, $z\in V$.
\end{theo}
\begin{proof}
Suppose that $R$ is a quadratic Leibniz conformal algebra. By
its definition, we set
\begin{eqnarray}
[x_{\lambda}y]=\partial(x\vdash y)+\lambda(x\vdash y+y\dashv x)+[x,y], ~~~\text{where $x$, $y\in V$,}
\end{eqnarray}
where $\vdash$, $\dashv$ and $[\cdot,\cdot]$ are three $\mathbb{C}$-bilinear maps
from $V\times V\rightarrow V$.

Next, we consider the Leibniz identity. For any $x$,
$y$ and $z\in V$, we get
\begin{eqnarray*}
 \ \ \ \ [x\,{}_{\lambda}\,{}[y\,{}_{\mu}\,{}z]]&=&[x\,{}_{\lambda}(\partial\,{}(y\vdash z)+\mu (y\vdash z+z\dashv y)+[y,z])]\nonumber\\
&=&(\lambda+\partial)[x\,{}_{\lambda}\,{}(y\vdash z)]+\mu[x\,{}_{\lambda}\,{}(y\vdash z+z\dashv y)]+[x\,{}_{\lambda}\,{}[y,z]]\\
&=&(\lambda+\partial)(\partial(x\vdash(y\vdash z))+\lambda (x\vdash(y\vdash z)+(y\vdash z)\dashv x)+[x, y\vdash z])\\
&&+\mu(\partial(x\vdash (y\vdash z+z\dashv y))+\lambda (x\vdash (y\vdash z+z\dashv y)+(y\vdash z+z\dashv y)\dashv x)\\
&&+[x, y\vdash z+z\dashv y])+\partial (x\vdash [y,z])+\lambda( x\vdash [y,z]+[y,z]\dashv x)+[x,[y,z]]\\
&=&\lambda \partial(2 x\vdash (y\vdash z)+(y\vdash z)\dashv x))+\lambda^{2}(x\vdash(y\vdash z)+(y\vdash z)\dashv x)\\
&&+\partial^{2}(x\vdash (y\vdash z))+\mu\partial(x\vdash (y\vdash z)+x\vdash (z\dashv y))\\
&&+\lambda\mu(x\vdash (y\vdash z)+x\vdash (z\dashv y)+(y\vdash z)\dashv x+(z\dashv y)\dashv x)\\
&&+\partial(x\vdash [y,z]+[x,y\vdash z])+\mu ([x,y\vdash z]+[x, z\dashv y])\\
&&+\lambda ([x,y\vdash z]+x\vdash [y,z]+[y,z]\dashv x)+[x,[y,z]].\\
\end{eqnarray*}
Similarly, we get
\begin{eqnarray*}
\ \ \ \ [y\,{}_{\mu}\,{}[x\,{}_{\lambda}\,{}z]]\\
&=&\mu \partial(2 y\vdash (x\vdash z)+(x\vdash z)\dashv y))+\mu^{2}(y\vdash(x\vdash z)+(x\vdash z)\dashv y)\\
&&+\partial^{2}(y\vdash (x\vdash z))+\lambda\partial(y\vdash (x\vdash z)+y\vdash (z\dashv x))\\
&&+\lambda\mu(y\vdash (x\vdash z)+y\vdash (z\dashv x)+(x\vdash z)\dashv y+(z\dashv x)\dashv y)\\
&&+\partial(y\vdash [x,z]+[y,x\vdash z])+\lambda ([y,x\vdash z]+[y, z\dashv x])\\
&&+\mu ([y,x\vdash z]+y\vdash [x,z]+[x,z]\dashv y)+[y,[x,z]].\\
\end{eqnarray*}

On the other hand, one can have
\begin{eqnarray*}
 \ \ \ \ [[x\,{}_{\lambda}\,{}y]\,{}_{\lambda+\mu}\,{}z]]\\
&=&\lambda^{2}((y\dashv x)\vdash z+z\dashv (y\dashv x))-\mu^{2}((x\vdash y)\vdash z+z\dashv (x\vdash y))\\
&&-\mu\partial (x\vdash y)\vdash z+\lambda\partial (y\dashv x)\vdash z\\
&&+\lambda\mu((y\dashv x)\vdash z+z\dashv (y\dashv x)-(x\vdash y)\vdash z-z\dashv (x\vdash y))\\
&&+\lambda([y\dashv x,z]+[x,y]\vdash z+z\dashv [x,y])+\partial [x,y]\vdash z\\
&&+\mu([x,y]\vdash z+z\dashv [x,y]-[x\vdash y,z])+[[x,y],z].
\end{eqnarray*}

By the Leibniz identity and comparing the coefficients of $\partial^2$, $\lambda\partial$, $\lambda^2$, $\lambda\mu$, $\mu\partial$,
$\mu^2$, $\partial$, $\lambda$,  $\mu$ and $\lambda^0\mu^0\partial^0$,  we get
\begin{gather}
\label{t10}x\vdash (y\vdash z)=y\vdash (x\vdash z),\\
\label{t11}2x\vdash (y\vdash z)+(y\vdash z)\dashv x=y\vdash (x\vdash z)+y\vdash (z\dashv x)+(y\dashv x)\vdash z,\\
\label{t12}x\vdash (y\vdash z)+(y\vdash z)\dashv x=(y\dashv x)\vdash z+z\dashv (y\dashv x),\\
x\vdash (y\vdash z)+x\vdash (z\dashv y)+(y\vdash z)\dashv x+(z\dashv y)\dashv x\nonumber\\
=y\vdash (x\vdash z)+y\vdash (z\dashv x)+(x\vdash z)\dashv y+(z\dashv x)\dashv y\nonumber\\
\label{t13}+(y\dashv x)\vdash z+z\dashv (y\dashv x)-(x\vdash y)\vdash z-z\dashv (x\vdash y),\\
\label{t14}x\vdash (y\vdash z)+x\vdash (z\dashv y)=2y\vdash (x\vdash z)+(x\vdash z)\dashv y-(x\vdash y)\vdash z,\\
\label{t15}y\vdash (x\vdash z)+(x\vdash z)\dashv y-(x\vdash y)\vdash z-z\dashv (x\vdash y)=0,\\
\label{t16}x\vdash[y,z]+[x,y\vdash z]=y\vdash [x,z]+[y, x\vdash z]+[x,y]\vdash z,\\
\label{t17}[x,y\vdash z]+x\vdash [y,z]+[y,z]\dashv x=[y,x\vdash z+z\dashv x]+[y\dashv x, z]+[x,y]\vdash z+z\dashv [x,y],\\
\label{t18}[x,y\vdash z+z\dashv y]=[y,x\vdash z]+y\vdash [x,z]+[x,z]\dashv y+[x,y]\vdash z+z\dashv [x,y]-[x\vdash y,z],\\
\label{t19}[x,[y,z]]=[y,[x,z]]+[[x,y],z].
\end{gather}

By (\ref{t19}), $(V, [\cdot,\cdot])$ is a Leibniz algebra. Note that (\ref{t10}) is just (\ref{t1}). Replacing $x$ and $y$ in (\ref{t14}) by $y$ and $x$ respectively and by (\ref{t11}), we can get (\ref{tt1}). By (\ref{t10}) and (\ref{tt1}), we can directly obtain (\ref{t3}) from (\ref{t11}). Moreover, by (\ref{t11}),  we can get (\ref{t2})
from (\ref{t12}). According to (\ref{t3}) and (\ref{tt1}),  (\ref{t4}) can be immediately deduced from (\ref{t15}).
By (\ref{t1}), (\ref{tt1}) and (\ref{t3}), from (\ref{t13}), we can obtain
\begin{eqnarray*}
0&=&(x\vdash (y\vdash z)-y\vdash (x\vdash z))+(x\vdash(z\dashv y)-(x\vdash z)\dashv y+(x\vdash y)\vdash z)\\
&&+(y\vdash z)\dashv x+(z\dashv y)\dashv x-y\vdash (z\dashv x)-(z\dashv x)\dashv y-(y\dashv x)\vdash z\\
&&-z\dashv(y\dashv x)+z\dashv (x\vdash y)\\
&=&(x\vdash (y\vdash z)+(y\vdash z)\dashv x-y\vdash (z\dashv x))-(y\dashv x)\vdash z+(z\dashv y)\dashv x \\
&&-(z\dashv x)\dashv y
+z\dashv(x\vdash y-y\vdash x)\\
&=&(y\vdash x-y\dashv x)\vdash z+(z\dashv y)\dashv x -(z\dashv x)\dashv y
+z\dashv(x\vdash y-y\vdash x)\\
&=& (z\dashv y)\dashv x -(z\dashv x)\dashv y
+z\dashv(x\vdash y-y\vdash x).
\end{eqnarray*}
This is just (\ref{t5}). Similarly, it is easy to show that (\ref{t7})-(\ref{t9}) are equivalent to (\ref{t16})-(\ref{t18}).

Now, the proof is finished.
\end{proof}
\begin{remark}\label{rem1}
By (\ref{t7}), it is easy to get that
\begin{eqnarray}
([x,y]+[y,x])\vdash z=0,
\end{eqnarray}
for any $x$, $y$ and $z\in V$.
\end{remark}

\begin{coro}\label{co1}
If for any $x$, $y\in V$, $x\vdash y=y\dashv x$ in the quadruple $(V, \vdash, \dashv, [\cdot,\cdot])$, then the conditions (\ref{t1})-(\ref{t9}) are equivalent to
the following equalities:
\begin{eqnarray}
\label{r1}&&x\vdash(y\vdash z)=y\vdash(x\vdash z)=(x\vdash y)\vdash z,\\
\label{r2}&&x\vdash [y,z]+[x,y\vdash z]=y\vdash[x,z]+[y,x\vdash z]+[x,y]\vdash z,\\
\label{r3}&&2x\vdash[y,z]=[y,x\vdash z]+[x\vdash y,z],\\
\label{r4}&&[x\vdash y-y\vdash x,z]=([x,y]+[y,x])\vdash z,
\end{eqnarray}
for $x$, $y$ and $z\in V$.
\end{coro}
\begin{proof}
It can be easily obtained from Theorem \ref{th2}.
\end{proof}
\begin{remark}
By Theorem \ref{th1} and Theorem \ref{th2},
if for any $x$, $y\in V$, $x\vdash y=y\dashv x$ in the quadruple $(V, \vdash, \dashv, [\cdot,\cdot])$ and $(V,\vdash )$ is not commutative, then the corresponding Leibniz conformal algebra is not a Lie conformal algebra. This provides
an useful method to construct Leibniz conformal algebras which are not
Lie conformal algebras.
\end{remark}
\begin{coro}\label{co2}
If for any $x$, $y\in V$,
$x\vdash y=x\dashv y$ in the quadruple $(V, \vdash, \dashv, [\cdot,\cdot])$, then the conditions (\ref{t1})-(\ref{t9}) are equivalent to that $(V,\vdash)$ is a Novikov algebra and
the following equalities hold:
\begin{eqnarray}
&&\label{k2}x\vdash [y,z]+[x,y\vdash z]=y\vdash [x,z]+[y,x\vdash z]+[x,y]\vdash z,\\
&&\label{k6}[y,x\vdash z]+[x\vdash z,y]=0,\\
&&\label{k7}x\vdash([y,z]+[z,y])=0,
\end{eqnarray}
for any $x$, $y$ and $z\in V$.
\end{coro}

\begin{proof}
It can be easily obtained from Theorem \ref{th2}.
\end{proof}
\begin{remark}
By Theorem \ref{th1}, a quadratic Lie conformal algebra is equivalent
to a Gel'fand-Dorfman algebra. If for any $x$, $y\in V$,
$x\vdash y=x\dashv y$ and the Leibniz algebra is
a Lie algebra in the quadruple $(V, \vdash, \dashv, [\cdot,\cdot])$, then the corresponding quadratic Leibniz conformal algebra is a Lie conformal algebra. (\ref{k6}) and (\ref{k7}) measure how far
a quadratic Leibniz conformal algebra is a Lie conformal algebra when $x\vdash y=x\dashv y$ in the quadruple $(V, \vdash, \dashv, [\cdot,\cdot])$ for any $x$, $y\in V$.
\end{remark}
\begin{coro}\label{co3}
If for any $x$, $y\in V$,
$x\vdash y=-y\dashv x$ in the quadruple $(V, \vdash, \dashv, [\cdot,\cdot])$,
then (\ref{t1})-(\ref{t9}) are equivalent to
the following equalities
\begin{eqnarray}
&&\label{m1}(x\vdash y)\vdash z=x\vdash(y\vdash z)=0,\\
&&\label{m2}x\vdash [y,z]+[x,y\vdash z]=y\vdash [x,z]+[y,x\vdash z]+[x,y]\vdash z,\\
&&\label{m3}[x\vdash y,z]=[y,x\vdash z],\\
&&\label{m4}[x\vdash y+y\vdash x,z]=-([x,y]+[y,x])\vdash z,
\end{eqnarray}
for any $x$, $y$ and $z\in V$. Moreover,
letting $[x,y]=a x\circ y-by\circ x$ for any $x$,
$y\in V$ and $a$, $b\in \mathbb{C}$, (\ref{m1})-(\ref{m4}) and the fact $[\cdot,\cdot]$ is a Leibniz bracket are equivalent to (\ref{m1}).
\end{coro}
\begin{proof}
It can be easily obtained from Theorem \ref{th2}.
\end{proof}
\begin{coro}\label{co4}
If $\vdash$ is trivial, i.e. for any $x$, $y\in V$,
$x\vdash y=0$ in the quadruple $(V, \vdash, \dashv, [\cdot,\cdot])$,
then (\ref{t1})-(\ref{t9}) are equivalent to
the following equalities
\begin{eqnarray}
&&\label{m5} x\dashv (y\dashv z)=0,\\
&&\label{m6}(x\dashv y)\dashv z=(x\dashv z)\dashv y,\\
&&\label{m7}[x\dashv y,z]=-z\dashv([x,y]+[y,x]),\\
&&\label{m8}[x,z]\dashv y=[x,z\dashv y]-z\dashv [x,y],
\end{eqnarray}
for any $x$, $y$, $z\in V$.
\end{coro}
\begin{proof}
It can be directly obtained from Theorem \ref{th2}.
\end{proof}

\begin{coro}\label{cco5}
If $\dashv$ is trivial, i.e. for any $x$, $y\in V$,
$x\dashv y=0$ in the quadruple $(V, \vdash, \dashv, [\cdot,\cdot])$,
then (\ref{t1})-(\ref{t9}) are equivalent to
the following equalities
\begin{eqnarray}
&&\label{m9} x\vdash (y\vdash z)=(x\vdash y)\vdash z=0,\\
&&\label{m10}x\vdash [y,z]=[x\vdash y,z]=0,\\
&&\label{m11}[x,y\vdash z]=[y,x\vdash z]+[x,y]\vdash z,
\end{eqnarray}
for any $x$, $y$, $z\in V$.
\end{coro}
\begin{proof}
It can be directly obtained from Theorem \ref{th2}.
\end{proof}

Next, we present an example to construct Leibniz conformal algebras according to Corollaries \ref{co3}, \ref{co4} and \ref{cco5}.
\begin{examp}
Let $V= \mathbb{C}a\oplus \mathbb{C}b$ be a 2-dimensional vector space. Define a bilinear operation $\star: V\times V\rightarrow V$ by
\begin{eqnarray}
a\star a=b,~~~a\star b=b\star a=b\star b=0.
\end{eqnarray}
Then $(V,\star)$ satisfies the equalities (\ref{m1}), (\ref{m9}) and (\ref{m5}), (\ref{m6}), if
we replace $\vdash$ and $\dashv$ by $\star$ respectively.

Therefore, according to Theorem \ref{th2} and Corollary \ref{co3},
we can get the following quadratic Leibniz conformal algebra
$R=\mathbb{C}[\partial]V=\mathbb{C}[\partial]a\oplus \mathbb{C}[\partial]b$
with the $\lambda$-product given by
\begin{eqnarray}
[a_\lambda a]=\partial b,~~~[a_\lambda b]=[b_\lambda a]=[b_\lambda b]=0.
\end{eqnarray}

By Theorem \ref{th2} and Corollary \ref{co4}, we obtain
a quadratic Leibniz conformal algebra $R=\mathbb{C}[\partial]V=\mathbb{C}[\partial]a\oplus \mathbb{C}[\partial]b$
with the $\lambda$-product
\begin{eqnarray}
[a_\lambda a]=\lambda b,~~~[a_\lambda b]=[b_\lambda a]=[b_\lambda b]=0.
\end{eqnarray}

Similarly, by Theorem \ref{th2} and Corollary \ref{cco5}, we can get
a quadratic Leibniz conformal algebra $R=\mathbb{C}[\partial]V=\mathbb{C}[\partial]a\oplus \mathbb{C}[\partial]b$
with the $\lambda$-product
\begin{eqnarray}
[a_\lambda a]=(\partial+\lambda)b,~~~[a_\lambda b]=[b_\lambda a]=[b_\lambda b]=0.
\end{eqnarray}
\end{examp}
\begin{remark}
In fact, the above examples can be generalized as follows.
Let $R=\mathbb{C}[\partial]V=\mathbb{C}[\partial]a\oplus \mathbb{C}[\partial]b$
with the $\lambda$-product as follows:
\begin{eqnarray}
[a_\lambda a]=f(\lambda,\partial)b,~~[a_\lambda b]=[b_\lambda a]=[b_\lambda b]=0.
\end{eqnarray}
Obviously, $R$ is a Leibniz conformal algebra. Moreover, it is easy to see
that $R$ is a Lie conformal algebra if and only if $f(\lambda,\partial)
=-f(-\lambda-\partial,\partial)$.

\end{remark}

Since in the case of Corollary \ref{co2}, $(V, \vdash, [\cdot,\cdot])$ is almost
a Gel'fand-Dorfman algebra, in the following, we mainly study the case in Corollary \ref{co1}.
\begin{defi}
Given a vector space $V$ with two bilinear operations $\vdash$ and $[\cdot,\cdot]$. We call $(V,\vdash)$ a \emph{Perm algebra} if $\vdash$ satisfies (\ref{r1}).
$(V,\vdash,[\cdot,\cdot])$ is called a \emph{Perm-Leibniz algebra}
if $(V,[\cdot,\cdot])$ is a Leibniz algebra and $\vdash$, $[\cdot,\cdot]$ satisfy (\ref{r1})-(\ref{r4}).
\end{defi}
\begin{remark}
The definition of Perm algebra can be referred to \cite{Ch1}.
\end{remark}
\begin{remark}
Obviously, for a Perm algebra $(V,\vdash)$,
$(V,\vdash,[\cdot,\cdot])$ is a Perm-Leibniz algebra with
$[a,b]=0$ for any $a$, $b\in V$. In this case, we say the quadratic Leibniz conformal algebra $R=\mathbb{C}[\partial]V$ corresponds to the Perm algebra $(V,\vdash)$.
\end{remark}
\begin{remark}
By Theorem \ref{th2}, the quadratic Leibniz conformal algebra corresponding
to the Perm-Leibniz algebra $(V,\vdash,[\cdot,\cdot])$ is
$R=\mathbb{C}[\partial]V$ with the $\lambda$-product given by
\begin{eqnarray}
[a_\lambda b]=\partial(a\vdash b)+2\lambda(a\vdash b)+[a,b],
\end{eqnarray}
for any $a$, $b\in V$.

Therefore, the coefficient algebra of $R=\mathbb{C}[\partial]V$ is an infinite-dimensional vector space $V\otimes \mathbb{C}[t,t^{-1}]$ with the Leibniz bracket given by
\begin{eqnarray}
[a\otimes t^m, b\otimes t^n]=(m-n)(a\vdash b)\otimes t^{m+n-1}
+[a,b]\otimes t^{m+n},~~~a,~~b\in V,~~m,~~n\in \mathbb{Z}.
\end{eqnarray}
\end{remark}

Next, we give a construction of Perm algebras from commutative associative algebras.

\begin{prop}\label{pr1}
Let $(A,\cdot)$ be a commutative associative algebra and
$P:A\rightarrow A$ be a linear operator satisfying
\begin{eqnarray}\label{ppp1}
P(P(x)\cdot y)=P(x)\cdot P(y), ~~x,~~y\in A.
\end{eqnarray}
Define
\begin{eqnarray}
x\circ y=P(x)\cdot y,~~~~x,~~y\in A.
\end{eqnarray}
Then $(A,\circ)$ is a Perm algebra.
\end{prop}
\begin{proof}
It is easy to check.
\end{proof}
\begin{remark}
$P:A\rightarrow A$ satisfying (\ref{ppp1}) is called \emph{an averaging operator} on $A$. More details on averaging operators can  refer to\cite{GP}. In addition, Proposition \ref{pr1} can be seen a particular case of the general statement on replicated algebras in \cite{gk}.
\end{remark}

Finally, we present some examples of Leibniz conformal algebras through the given constructions.

\begin{examp}\label{exs1}
Note that the classification of 2-dimensional Novikov algebras up to isomorphism has been given in \cite{BM}.
Set $V=\mathbb{C}e_1\oplus \mathbb{C}e_2$. There are the following cases up to isomorphism.\\
(A) Communicative associative algebras:\\
(A1) \begin{eqnarray}
e_1\vdash e_1=e_1\vdash e_2=e_2\vdash e_1=e_2\vdash e_2=0;
\end{eqnarray}
(A2)\begin{eqnarray}
e_1\vdash e_1=e_2,~~e_1\vdash e_2=e_2\vdash e_1=e_2\vdash e_2=0;
\end{eqnarray}
(A3)\begin{eqnarray}
e_1\vdash e_1=e_1,~~e_1\vdash e_2=e_2\vdash e_1=0,~~e_2\vdash e_2=e_2;
\end{eqnarray}
(A4)\begin{eqnarray}
e_1\vdash e_1=e_1,~~e_1\vdash e_2=e_2\vdash e_1=e_2\vdash e_2=0;
\end{eqnarray}
(A5)\begin{eqnarray}
e_1\vdash e_1=e_1,~~e_1\vdash e_2=e_2\vdash e_1=e_2,~~e_2\vdash e_2=0;
\end{eqnarray}
(B) Associative but non-communicative Novikov algebra:
\begin{eqnarray}
e_1\vdash e_1=e_1\vdash e_2=0,~~e_2\vdash e_1=e_1,~~e_2\vdash e_2=e_2;
\end{eqnarray}
(C) Nonassociative Novikov algebras:\\
(C1)\begin{eqnarray}
e_1\vdash e_1=e_2\vdash e_1=e_2\vdash e_2= 0,~~e_1\vdash e_2=-e_1;
\end{eqnarray}
(C2)\begin{eqnarray}
e_1\vdash e_1=e_1\vdash e_2=0, ~~e_2\vdash e_1=e_1,~~e_2\vdash e_2=e_1+e_2;
\end{eqnarray}
(C3(l))\begin{eqnarray}
e_1\vdash e_1=0,~~e_1\vdash e_2=le_1, ~~e_2\vdash e_1=e_1,~~e_2\vdash e_2=e_2,~~l\neq 0, 1.
\end{eqnarray}
It is easy to see that if the Novikov algebras are isomorphic, then the corresponding quadratic Leibniz conformal algebras are isomorphic.

Note that (A1)-(A5) and (B) are Perm algebras. From a Perm algebra $(V,\vdash)$, we can naturally obtain a quadratic Leibniz conformal algebra $R=\mathbb{C}[\partial]V$ with the $\lambda$-product $[a_\lambda b]=\partial(a\vdash b)+2\lambda(a\vdash b)$ for any $a$, $b\in V$.
By Corollary \ref{co1}, the quadratic Leibniz conformal algebras corresponding to (A1)-(A5) are Lie conformal algebras and the quadratic Leibniz conformal algebra corresponding to (B) is not a Lie conformal algebra. We denote the Lie conformal algebra corresponding to (Ai) by $LC(Ai)$ for $i\in \{1,\cdots,5\}$ and the Leibniz conformal algebra corresponding to (B) by $LBC(B)$.

In addition, for any Novikov algebra $(V,\vdash)$, we can naturally obtain a quadratic Leibniz conformal algebra $R=\mathbb{C}[\partial]V$ with the $\lambda$-product $[a_\lambda b]=\partial(a\vdash b)+\lambda(a\vdash b+b\vdash a)$ for any $a$, $b\in V$.
By Corollary \ref{co2}, the quadratic Leibniz conformal algebras corresponding to (A1)-(A5), (B), (C1)-(C3) are Lie conformal algebras. We also denote the Lie conformal algebras corresponding to (B), (C1)-(C3(l)) by $LC(B)$  and $LC(C1)$-$LC(C3(l))$ respectively.

Here, we only present the definitions of $LBC(B)$ and $LC(B)$ in details.

$LBC(B)=\mathbb{C}[\partial]e_1\oplus \mathbb{C}[\partial]e_2$ is a Leibniz conformal algebra with the following $\lambda$-product:
\begin{eqnarray}
[{e_1}_\lambda e_1]=[{e_1}_\lambda e_2]=0,~~[{e_2}_\lambda e_1]=(\partial+2\lambda)e_1,~~[{e_2}_\lambda e_2]=(\partial+2\lambda)e_2.
\end{eqnarray}
The coefficient algebra of $LBC(B)$ is an infinite-dimensional vector space $V\otimes \mathbb{C}[t,t^{-1}]$ with the Leibniz bracket given by
\begin{eqnarray*}
&&[e_1\otimes t^m, e_1\otimes t^n]=[e_1\otimes t^m, e_2\otimes t^n]=0,\\
&&[e_2\otimes t^m, e_1\otimes t^n]=(m-n)e_1\otimes t^{m+n-1},\\
&&[e_2\otimes t^m, e_2\otimes t^n]=(m-n)e_2\otimes t^{m+n-1},
\end{eqnarray*}
for any $m$, $n\in \mathbb{Z}$.

$LC(B)=\mathbb{C}[\partial]e_1\oplus \mathbb{C}[\partial]e_2$ is a Lie conformal algebra with the following $\lambda$-bracket:
\begin{eqnarray}
[{e_1}_\lambda e_1]=0,~~[{e_1}_\lambda e_2]=\lambda e_1,~~[{e_2}_\lambda e_1]=(\partial+\lambda)e_1,~~[{e_2}_\lambda e_2]=(\partial+2\lambda)e_2.
\end{eqnarray}
The coefficient algebra of $LC(B)$ is an infinite-dimensional vector space $V\otimes \mathbb{C}[t,t^{-1}]$ with the Lie bracket given by
\begin{eqnarray*}
&&[e_1\otimes t^m, e_1\otimes t^n]=0,\\
&&[e_1\otimes t^m, e_2\otimes t^n]=me_1\otimes t^{m+n-1},\\
&&[e_2\otimes t^m, e_2\otimes t^n]=(m-n)e_2\otimes t^{m+n-1},
\end{eqnarray*}
for any $m$, $n\in \mathbb{Z}$.
\end{examp}
\section{Central extensions of quadratic Leibniz conformal algebras}

In this section, we will study central extensions of quadratic Leibniz conformal algebras.

An extension of a Leibniz conformal algebra $R$ by an abelian Leibniz conformal algebra $C$ is a short exact sequence of Leibniz conformal algebras
\begin{eqnarray*}
0\rightarrow C\rightarrow \widehat{R}\rightarrow R\rightarrow 0.
\end{eqnarray*}
$\widehat{R}$ is called an \emph{extension} of $R$ by $C$ in this case. This extension is called \emph{central} if $\partial C=0$ and
$[C_\lambda \widehat{R}]=0$. Moreover, if $R$ is a Lie conformal algebra, then we call this extension
a \emph{Leibniz central extension} of this Lie conformal algebra.

In the following, we investigate the central extensions $\widehat{R}$ of $R$ by a one-dimensional center $\mathbb{C}\mathfrak{c}$. This implies that $\widehat{R}=R\oplus \mathbb{C}\mathfrak{c}$, and
\begin{eqnarray*}
[a_\lambda b]_{\widehat{R}}=[a_\lambda b]_R+\alpha_\lambda(a,b)\mathfrak{c}, \text{for~~all~~$a$, $b\in R$,}
\end{eqnarray*}
where $\alpha_\lambda(\cdot,\cdot): R\times R\rightarrow \mathbb{C}[\lambda]$ is a $\mathbb{C}$-bilinear map.
By the axioms of Leibniz conformal algebra, $\alpha_\lambda(\cdot,\cdot)$ should satisfy the following properties (for all $a$, $b$, $c\in R$) :
\begin{eqnarray}
&&\label{f1}\alpha_\lambda(\partial a,b)=-\lambda \alpha_\lambda(a,b)=-\alpha_\lambda(a,\partial b),\\
&&\label{f2}\alpha_\lambda(a,[b_\mu c])=\alpha_{\lambda+\mu}([a_\lambda b],c)+\alpha_{\mu}(b,[a_\lambda c]).
\end{eqnarray}
By the cohomology theory of Leibniz conformal algebras in \cite{BKV}, $\alpha_\lambda(\cdot,\cdot)$ is a 2-cocycle in the reduced complex of $R$ with values in the trivial $R$-module $\mathbb{C}$ and the central extensions of $R$ by $\mathbb{C}\mathfrak{c}$ up to equivalence can be parameterized by $H^2(R,\mathbb{C}\mathfrak{c})$. According to the definition of 2-coboundary, 2-cocycles $\alpha_\lambda(\cdot,\cdot)$ and $\alpha_\lambda^{'}(\cdot,\cdot)$ are equivalent if and only if there exists a $\mathbb{C}[\partial]$-module homomorphism $\varphi: R\rightarrow \mathbb{C}$ ($\mathbb{C}$ is a trivial $\mathbb{C}[\partial]$-module) such that $\alpha_\lambda(a,b)=\alpha_\lambda^{'}(a,b)+\varphi([a_\lambda b])$ for all
$a$, $b\in R$. Equivalent 2-cocycles define equivalent central extensions.

For the quadratic Leibniz conformal algebra corresponding to $(V, \vdash, \dashv, [\cdot,\cdot])$, we can directly obtain the following conformal cocycle condition from (\ref{f2}):
\begin{gather}
(\lambda+\mu) \alpha_\lambda(a,b\vdash c)+\mu \alpha_\lambda(a,c\dashv b)+\alpha_\lambda(a,[b,c])\nonumber\\
=-\mu\alpha_{\lambda+\mu}(a\vdash b,c)+\lambda\alpha_{\lambda+\mu}(b\dashv a,c)+\alpha_{\lambda+\mu}([a,b],c)\nonumber\\
\label{cy1}+(\lambda+\mu) \alpha_\mu(b,a\vdash c)+\lambda \alpha_\mu(b,c\dashv a)+\alpha_\mu(b,[a,c]).
\end{gather}

\begin{theo}\label{th3}
Let $V$ be a vector space and $\widehat{R}=R\oplus\mathbb{C}\mathfrak{c}$ be a central extension of quadratic Leibniz conformal algebra $R=\mathbb{C}[\partial]V$  corresponding to the Perm-Leibniz algebra $(V,\vdash, [\cdot,\cdot])$ by a one-dimensional center $\mathbb{C}\mathfrak{c}$.
Set the $\lambda$-product of $\widehat{R}$ be
\begin{eqnarray}
\label{f3}\widetilde{[a_\lambda b]}=\partial(a\vdash b)+2\lambda(a\vdash b)+[a,b]+\alpha_\lambda(a,b)\mathfrak{c},
\end{eqnarray}
where $a$, $b\in V$  and $\alpha_\lambda(a,b)\in \mathbb{C}[\lambda]$.
Assume that $\alpha_\lambda(a,b)=\sum_{i=0}^n\lambda^i\alpha_i(a,b)$ for any $a$, $b\in V$, where $\alpha_i(\cdot,\cdot): V\times V\rightarrow \mathbb{C}$ are bilinear forms and there exist some $a$, $b \in V$ such that
$\alpha_n(a,b)\neq 0$. Then we obtain (for any $a$, $b$, $c\in V$):\\
(1) If $n>3$,  then $\alpha_n(a\vdash b, c)=0$ ;\\
(2) If $n\leq 3$, then $\alpha_\lambda(a,b)=\sum_{i=0}^3\lambda^i\alpha_i(a,b)$ and $\alpha_\lambda(\cdot,\cdot)$ satisfies the following conformal cocycle condition:
\begin{gather}
(\lambda+2\mu) \alpha_\lambda(a,b\vdash c)+\alpha_\lambda(a,[b,c])\nonumber\\
\label{cy2}=(\lambda-\mu)\alpha_{\lambda+\mu}(a\vdash b,c)+\alpha_{\lambda+\mu}([a,b],c)+(2\lambda+\mu) \alpha_\mu(b,a\vdash c)+\alpha_\mu(b,[a,c]),
\end{gather}
for any $a$, $b$, $c\in V$;\\
(3) Such two 2-cocycles $\alpha_\lambda(\cdot,\cdot)$ and $\alpha_\lambda^{'}(\cdot,\cdot)$ are equivalent if and only if there exists a linear map
$\varphi: V\rightarrow \mathbb{C}$ such that
\begin{eqnarray}
\alpha_\lambda(a,b)=\alpha_\lambda^{'}(a,b)+2\lambda \varphi(a\vdash b)+\varphi([a,b]), ~~~~\text{for all $a$, $b\in V$}.
\end{eqnarray}

\end{theo}
\begin{proof}
Note that (2) and (3) can be directly obtained from (\ref{cy1}) and the definition of equivalence of 2-cocycles. Therefore, we only need to prove (1).

If $n>3$, by the assumption of $\alpha_\lambda(a,b)=\sum_{i=0}^n\lambda^i\alpha_i(a,b)$ and comparing the coefficients of $\lambda^2\mu^{n-1}$ in (\ref{cy2}), we get
\begin{eqnarray}
(n-C_n^2)\alpha_n(a\vdash b,c)=0.
\end{eqnarray}
Therefore, $\alpha_n(a\vdash b,c)=0$.

This completes the proof.
\end{proof}
\begin{coro}\label{co5}
Let $R=\mathbb{C}[\partial]V$ be a finite quadratic Leibniz conformal algebra
corresponding to the Perm-Leibniz algebra $(V,\vdash,[\cdot,\cdot])$ where $V=V\vdash V$. Set $\widehat{R}=R\oplus
\mathbb{C}\mathfrak{c}$ be a central extension of $(R,[\cdot_\lambda \cdot])$ with the $\lambda$-product
given by (\ref{f3}).
Then for any $a$, $b\in V$, we obtain that $\alpha_\lambda(a,b)
=\sum_{i=0}^3\lambda^i\alpha_i(a,b)$ where $\alpha_i(\cdot,\cdot): V\times V\rightarrow \mathbb{C}$ are bilinear forms and $\alpha_\lambda(\cdot,\cdot)$ satisfies (\ref{cy2}).
\end{coro}
\begin{proof}
Since $R=\mathbb{C}[\partial]V$ is a finitely generated and free $\mathbb{C}[\partial]$-module, we can assume that $\alpha_\lambda(a,b)=\sum_{i=0}^n\lambda^i\alpha_i(a,b)$ for any $a$, $b\in V$ and some non-negative integer $n$.
Since $V=V\vdash V$, for any $z\in V$, there exist some $m$ and $x_i$, $y_i\in V$ such that $z=\sum_{i=0}^mx_i\circ y_i$.
Therefore, by Theorem \ref{th3}, we get that if $n>3$, for any $z\in V$ and $c\in V$,
$\alpha_n(z,c)=\sum_{i=0}^m\alpha_n(x_i\circ y_i,c)=0$. Hence, for any $a$, $b\in V$, we obtain $\alpha_\lambda(a,b)
=\sum_{i=0}^3\lambda^i\alpha_i(a,b)$. By Theorem \ref{th3}, we can obtain this corollary.
\end{proof}

\begin{remark}
It should be pointed out that there are some natural conditions to make the Perm algebra $(V,\vdash)$ satisfy the assumption in Corollary \ref{co5}. For example, $(V,\vdash)$ is simple or $(V,\vdash)$ has a left unit or a right unit.
\end{remark}

\begin{coro}\label{co6}
Let $R=\mathbb{C}[\partial]V$ be a quadratic Leibniz conformal algebra
corresponding to the Perm algebra $(V,\vdash)$ where $V=V\vdash V$. Set $\widehat{R}=R\oplus
\mathbb{C}\mathfrak{c}$ be a central extension of $(R,[\cdot_\lambda \cdot])$ with the $\lambda$-product
given as follows:
\begin{eqnarray}
\label{cy3}\widetilde{[a_\lambda b]}=\partial(a\vdash b)+2\lambda(a\vdash b)+\alpha_\lambda(a,b)\mathfrak{c}, ~~~~~~~~~~~~~~~~a,~~b\in V.
\end{eqnarray}
Then for any $a$, $b\in V$, we obtain that $\alpha_\lambda(a,b)
=\alpha_0(a,b)+\lambda\alpha_1(a,b)+\lambda^3\alpha_3(a,b)$ and the bilinear forms  $\alpha_0(\cdot,\cdot)$,
$\alpha_1(\cdot,\cdot)$ and $\alpha_3(\cdot,\cdot)$ on $V$ satisfy the following
equalities:
\begin{eqnarray}
&&\label{b1}\alpha_3(a,b\vdash c)=\alpha_3(a\vdash b,c)=\alpha_3(b, a\vdash c),\\
&&\label{b2}\alpha_1(a,b\vdash c)
=\alpha_1(a\vdash b,c)=\alpha_1(b, a\vdash c),\\
&&\label{b9}\alpha_0(a,b\vdash c)
=-\alpha_0(a\vdash b,c)=\alpha_0(b, a\vdash c),
\end{eqnarray}
for any $a$, $b$, $c\in V$.
Moreover,  such two 2-cocycles $\alpha_\lambda(\cdot,\cdot)$ and $\alpha_\lambda^{'}(\cdot,\cdot)$ are equivalent if and only if $\alpha_i(a,b)=\alpha_i^{'}(a,b)$ for $i=0$, and $i=3$, and there exists a linear map
$\varphi: V\rightarrow \mathbb{C}$ such that
\begin{eqnarray}
\alpha_1(a,b)=\alpha_1^{'}(a,b)+2\varphi(a\vdash b),
\end{eqnarray}
for all $a$, $b\in V$.
\end{coro}
\begin{proof}
For any $a$, $b\in V$, set $\alpha_\lambda(a,b)=\sum_{i=0}^{n_{a,b}}\lambda^i\alpha_i(a,b)$ where
$\alpha_i(\cdot,\cdot)$ are bilinear forms on $V$ and $n_{a,b}$ is a non-negative integer depending on $a$ and $b$.
The conformal cocycle condition becomes
\begin{gather}
\label{o1}(\lambda+2\mu) \alpha_\lambda(a,b\vdash c)
=(\lambda-\mu)\alpha_{\lambda+\mu}(a\vdash b,c)+(2\lambda+\mu) \alpha_\mu(b,a\vdash c).
\end{gather}
For fixed $a$, $b$, $c$, there are only finite elements of $V$ appearing in $\alpha_\lambda(\cdot,\cdot)$ in (\ref{o1}). Therefore, we may assume the degrees of all $\alpha_\lambda(\cdot,\cdot)$ in (\ref{o1}) are smaller than some non-negative integer. So, we set $\alpha_\lambda(a,b\vdash c)
=\sum_{i=0}^n\lambda^i\alpha_i(a, b\vdash c)$, $\alpha_{\lambda+\mu}(a\vdash b,c)=\sum_{i=0}^n(\lambda+\mu)^i\alpha_i(a\vdash b,c)$ and
$\alpha_\mu(b,a\vdash c)=\sum_{i=0}^n\mu^i\alpha_i(b,a\vdash c)$.

If $n>3$, by comparing the coefficients of $\lambda^2\mu^{n-1}$ in (\ref{o1}), we get
\begin{eqnarray}
(n-C_n^2)\alpha_n(a\vdash b,c)=0.
\end{eqnarray}
Therefore, $\alpha_n(a\vdash b,c)=0$. Repeating this process, we can get
$\alpha_m(a\vdash b,c)=0$ for all $n\geq m> 3$.

Based on the above, for any $a$, $b$, $c\in V$, we get $\alpha_m(a\vdash b,c)=0$ for all $m>3$.
Since for any $z\in V$, there exist some $m$ and $x_i$, $y_i\in V$ such that $z=\sum_{i=0}^m x_i\circ y_i$, we get
$\alpha_m(z,c)=0$ for all $m>3$. Hence, $\alpha_\lambda(a,b)=\sum_{i=0}^3 \lambda^i\alpha_i(a,b)$. Taking it into (\ref{o1}) and by comparing the coefficients of $\lambda^4$, $\lambda^3\mu$, $\lambda^2\mu^2$, $\lambda\mu^3$, $\mu^4$,
$\lambda^3$, $\lambda^2\mu$, $\lambda\mu^2$, $\mu^3$, $\lambda^2$, $\lambda\mu$, $\mu^2$, $\lambda$, $\mu$ and and $\lambda^0\mu^0$, it is easy to check that the conformal cocycle condition (\ref{o1}) is equivalent to that $\alpha_2(\cdot,\cdot)=0$ and $\alpha_3(\cdot,\cdot)$, $\alpha_1(\cdot,\cdot)$, $\alpha_0(\cdot,\cdot)$ satisfy (\ref{b1})-(\ref{b9}).
\end{proof}

\begin{remark}
Note that this corollary also holds for the quadratic Leibniz conformal algebras corresponding to Perm algebras which are not finite.
\end{remark}

\begin{coro}\label{co7}
Let $(V,\vdash)$ be a Perm algebra and $R=\mathbb{C}[\partial]V$ be the corresponding quadratic Leibniz conformal algebra. Set $\alpha_i(\cdot,\cdot)$ $(i=0, 1, 3)$
be bilinear forms on $V$  satisfying (\ref{b1})-(\ref{b9}). Define bilinear forms $\varphi_i: \text{Coeff}(R)\times \text{Coeff}(R) \rightarrow \mathbb{C}$ as follows:
\begin{eqnarray}
&&\varphi_0(a\otimes t^m, b\otimes t^n)=\alpha_0(a,b)\delta_{m+n+1,0},\\
&&\varphi_1(a\otimes t^m, b\otimes t^n)=m\alpha_1(a,b)\delta_{m+n,0},\\
&&\varphi_3(a\otimes t^m, b\otimes t^n)=m(m-1)(m-2)\alpha_3(a,b)\delta_{m+n-2,0},
\end{eqnarray}
for $a$, $b\in V$, $m$, $n\in \mathbb{Z}$. Then $\varphi_0$, $\varphi_1$, $\varphi_3$ are 2-cocycles of Leibniz algebra $\text{Coeff}(R)$.
\end{coro}
\begin{proof}
Suppose that $\widehat{R}$ is the central extension defined by (\ref{cy3}). Then the coefficient algebra of $\widehat{R}$ is $\text{Coeff}(\widehat{R})=\text{Coeff}(R)\oplus \mathbb{C}\mathfrak{c}_{-1}$ with the following Leibniz bracket:
\begin{eqnarray*}[a_m, b_n]&=&(m-n)(a\vdash b)_{m+n-1}+\alpha_0(a,b)\delta_{m+n+1,0}{\mathfrak{c}}_{-1}\\
&&+m\alpha_1(a,b)\delta_{m+n,0}{\mathfrak{c}}_{-1}+m(m-1)(m-2)\alpha_3(a,b)\delta_{m+n-2,0}{\mathfrak{c}}_{-1},\end{eqnarray*}
for any $a$, $b\in V$, $m$, $n\in \mathbb{Z}$, and $\mathfrak{c}_{-1}$ is in the center of $\text{Coeff}(\widehat{R})$. Obviously,  $\text{Coeff}(\widehat{R})$ is the central extension of $\text{Coeff}(R)$ by a one-dimensional center $\mathbb{C}\mathfrak{c}_{-1}$. Since $\alpha_0(\cdot,\cdot)$, $\alpha_1(\cdot,\cdot)$ and
$\alpha_3(\cdot,\cdot)$ donot depend on each other in (\ref{b1})-(\ref{b9}),
we naturally obtain 3 kinds of 2-cocycles of $\text{Coeff}(R)$.
\end{proof}

\begin{examp}
In this example, we give a characterization of the central extensions
of $LC(A1)$- $LC(A5)$ and $LBC(B)$ given in Example \ref{exs1} by a one-dimensional center $\mathbb{C}\mathfrak{c}$ up to equivalence.

Obviously, the Perm algebras $(A3)$, $(A5)$ and $(B)$ satisfy the condition
$V=V\vdash V$ in Corollary \ref{co6}. For investigating the central extensions
of $LC(A3)$, $LC(A5)$ and $LBC(B)$ by a one-dimensional center $\mathbb{C}\mathfrak{c}$, we only need to compute the bilinear forms
$\alpha_0(\cdot,\cdot)$, $\alpha_1(\cdot,\cdot)$, $\alpha_3(\cdot,\cdot)$ on $V$ satisfying (\ref{b1})-(\ref{b9}).

First, we consider the case of $LBC(B)$. By the definition of (B), it is easy to get that $\alpha_3(e_1,e_1)=\alpha_3(e_1,e_2)=\alpha_1(e_1,e_1)=\alpha_1(e_1,e_2)=0$,
$\alpha_3(e_2,e_1)=\beta_1$, $\alpha_3(e_2,e_2)=\beta_2$, $\alpha_1(e_2,e_1)=\gamma_1$, $\alpha_1(e_2,e_2)=\gamma_2$ and $\alpha_0(e_1,e_1)=\alpha_0(e_1,e_2)=\alpha_0(e_2,e_1)=\alpha_0(e_2,e_2)=0$ where $\beta_1$, $\beta_2$, $\gamma_1$, $\gamma_2\in \mathbb{C}$. By choosing the linear map $\varphi:V\rightarrow\mathbb{C}$ in Corollary \ref{co6} defined by
$\varphi(e_1)=\frac{\gamma_1}{2}$ and $\varphi(e_2)=\frac{\gamma_2}{2}$, we can make $\alpha_1(\cdot,\cdot)$ be zero up to equivalence. Therefore, by Corollary \ref{co6}, all equivalence classes of central extensions of
$LBC(B)$ by a one-dimensional center $\mathbb{C}\mathfrak{c}$ are $\widehat{LBC(B)}(\beta_1,\beta_2)$  with the $\lambda$-product as follows:
\begin{eqnarray}
&&[{e_1}_\lambda e_1]=[{e_1}_\lambda e_2]=0,~~[{e_2}_\lambda e_1]=(\partial+2\lambda)e_1+\beta_1\lambda^3\mathfrak{c},\\
&&[{e_2}_\lambda e_2]=(\partial+2\lambda)e_2+\beta_2\lambda^3\mathfrak{c},
\end{eqnarray}
for all $\beta_1$, $\beta_2\in \mathbb{C}$. Obviously, if $(\beta_1,\beta_2)\neq (\beta_1^{'},\beta_2^{'})$, then $\widehat{LBC(B)}(\beta_1,\beta_2)$ is not equivalent to $\widehat{LBC(B)}(\beta_1^{'},\beta_2^{'})$. Therefore, $\text{dim}~H^2(LBC(B),\mathbb{C}\mathfrak{c})=2$.

With a similar discussion as above, we can get the central extensions
of $LC(A3)$ and $LC(A5)$ by a one-dimensional center $\mathbb{C}\mathfrak{c}$ up to equivalence as follows.

All equivalence classes of central extensions of
$LC(A3)$ by a one-dimensional center $\mathbb{C}\mathfrak{c}$ are $\widehat{LC(A3)}(\beta_1,\beta_2)$  with the $\lambda$-product as follows:
\begin{eqnarray}
&&[{e_1}_\lambda e_1]=(\partial+2\lambda)e_1+\beta_1\lambda^3\mathfrak{c},~~[{e_1}_\lambda e_2]=[{e_2}_\lambda e_1]=0,\\
&&[{e_2}_\lambda e_2]=(\partial+2\lambda)e_2+\beta_2\lambda^3\mathfrak{c},
\end{eqnarray}
for all $\beta_1$, $\beta_2\in \mathbb{C}$. Obviously, if $(\beta_1,\beta_2)\neq (\beta_1^{'},\beta_2^{'})$, then $\widehat{LC(A3)}(\beta_1,\beta_2)$ is not equivalent to $\widehat{LC(A3)}(\beta_1^{'},\beta_2^{'})$. Therefore, in this case, $\text{dim}~H^2(LC(A3),\mathbb{C}\mathfrak{c})=2$.

All equivalence classes of central extensions of
$LC(A5)$ by a one-dimensional center $\mathbb{C}\mathfrak{c}$ are $\widehat{LC(A5)}(\beta_1,\beta_2)$  with the $\lambda$-product as follows:
\begin{eqnarray}
&&[{e_1}_\lambda e_1]=(\partial+2\lambda)e_1+\beta_1\lambda^3\mathfrak{c},~~[{e_2}_\lambda e_2]=0,\\
&&[{e_1}_\lambda e_2]=[{e_2}_\lambda e_1]=(\partial+2\lambda)e_2+\beta_2\lambda^3\mathfrak{c},
\end{eqnarray}
for all $\beta_1$, $\beta_2\in \mathbb{C}$. Obviously, if $(\beta_1,\beta_2)\neq (\beta_1^{'},\beta_2^{'})$, then $\widehat{LC(A5)}(\beta_1,\beta_2)$ is not equivalent to $\widehat{LC(A5)}(\beta_1^{'},\beta_2^{'})$. Therefore, in this case, $\text{dim}~H^2(LC(A5),\mathbb{C}\mathfrak{c})=2$.

Finally, we determine the central extensions of $LC(A1)$, $LC(A2)$ and $LC(A4)$ by a one-dimensional center $\mathbb{C}\mathfrak{c}$ up to equivalence. By the conformal cocycle condition (\ref{o1}), it is easy to see that
$\text{dim}~H^2(LC(A1),\mathbb{C}\mathfrak{c})=\infty$, $\text{dim}~H^2(LC(A2),\mathbb{C}\mathfrak{c})=\infty$ and $\text{dim}~H^2(LC(A4),\mathbb{C}\mathfrak{c})=\infty$.
\end{examp}

Next, we study the Leibniz central extensions of quadratic Lie conformal algebras.
\begin{theo}\label{th4}
Let $V$ be a vector space and $\widehat{R}=R\oplus\mathbb{C}\mathfrak{c}$ be a Leibniz central extension of quadratic Lie conformal algebra $R=\mathbb{C}[\partial]V$  corresponding to the Gel'fand-Dorfman algebra $(V,\vdash, [\cdot,\cdot])$ by a one-dimensional center $\mathbb{C}\mathfrak{c}$.
Set the $\lambda$-bracket of $\widehat{R}$ be
\begin{eqnarray}
\label{o4}\widetilde{[a_\lambda b]}=\partial(a\vdash b)+\lambda(a\vdash b+b\vdash a)+[a,b]+\alpha_\lambda(a,b)\mathfrak{c},
\end{eqnarray}
where $a$, $b\in V$ and $\alpha_\lambda(a,b)\in \mathbb{C}[\lambda]$. Assume that $\alpha_\lambda(a,b)=\sum_{i=0}^n\lambda^i\alpha_i(a,b)$ for any $a$, $b\in V$, where all $\alpha_i(\cdot,\cdot): V\times V\rightarrow \mathbb{C}$ are bilinear forms and there exist some $a$, $b \in V$ such that
$\alpha_n(a,b)\neq 0$. Then we obtain (for any $a$, $b$, $c\in V$):\\
(1) If $n>3$,  then $\alpha_n(a\vdash b, c)=0$ ;\\
(2) If $n\leq 3$,  then $\alpha_\lambda(a,b)=\sum_{i=0}^3\lambda^i\alpha_i(a,b)$ and $\alpha_\lambda(\cdot,\cdot)$ satisfies the following conformal cocycle condition:
\begin{gather}
(\lambda+\mu) \alpha_\lambda(a,b\vdash c)+\mu \alpha_\lambda(a,c\vdash b)+\alpha_\lambda(a,[b,c])\nonumber\\
=-\mu\alpha_{\lambda+\mu}(a\vdash b,c)+\lambda\alpha_{\lambda+\mu}(b\vdash a,c)+\alpha_{\lambda+\mu}([a,b],c)\nonumber\\
\label{cy4}+(\lambda+\mu) \alpha_\mu(b,a\vdash c)+\lambda \alpha_\mu(b,c\vdash a)+\alpha_\mu(b,[a,c]).
\end{gather}
for any $a$, $b$, $c\in V$;\\
(3) Such two 2-cocycles $\alpha_\lambda(\cdot,\cdot)$ and $\alpha_\lambda^{'}(\cdot,\cdot)$ are equivalent if and only if there exists a linear map
$\varphi: V\rightarrow \mathbb{C}$ such that
\begin{eqnarray}
\alpha_\lambda(a,b)=\alpha_\lambda^{'}(a,b)+\lambda \varphi(a\vdash b+b\vdash a)+\varphi([a,b]), ~~~~\text{for all $a$, $b\in V$}.
\end{eqnarray}
\end{theo}
\begin{proof}
The proof is similar to that in Theorem \ref{th3} and can be referred to that in Theorem 3.1 in \cite{H1}.
\end{proof}

\begin{coro}
Let $V$ be a vector space and $R=\mathbb{C}[\partial]V$ be a finite quadratic Lie conformal algebra
corresponding to the Gel'fand-Dorfman algebra $(V,\vdash,[\cdot,\cdot])$ where $V=V\vdash V$. Set $\widehat{R}=R\oplus
\mathbb{C}\mathfrak{c}$ be a Leibniz central extension of $(R,[\cdot_\lambda \cdot])$ with the $\lambda$-product
given by (\ref{o4}).
Then for any $a$, $b\in V$, we obtain that $\alpha_\lambda(a,b)
=\sum_{i=0}^3\lambda^i\alpha_i(a,b)$ where $\alpha_i(\cdot,\cdot): V\times V\rightarrow \mathbb{C}$ are bilinear forms and $\alpha_\lambda(\cdot,\cdot)$ satisfies (\ref{cy4}).
\end{coro}
\begin{proof}
It can be obtained from Theorem \ref{th4} with a similar proof as that in Corollary \ref{co5}.
\end{proof}
\begin{coro}\label{cp1}
Let $R=\mathbb{C}[\partial]V$ be a quadratic Lie conformal algebra
corresponding to the Novikov algebra $(V,\vdash)$ where $V=V\vdash V$. Set $\widehat{R}=R\oplus
\mathbb{C}\mathfrak{c}$ be a Leibniz central extension of $(R,[\cdot_\lambda \cdot])$ with the $\lambda$-product
given as follows:
\begin{eqnarray}
\widetilde{[a_\lambda b]}=\partial(a\vdash b)+\lambda(a\vdash b+b\vdash a)+\alpha_\lambda(a,b)\mathfrak{c}, ~~~~~~~~~~~~~~~~
\end{eqnarray}
 for any $a$, $b\in V$.
Then for any $a$, $b\in V$, we obtain that $\alpha_\lambda(a,b)
=\sum_{i=0}^3\lambda^i\alpha_i(a,b)$ and the bilinear forms $\alpha_i(\cdot,\cdot)$ satisfy  the following equalities for any $a$, $b$, $c\in V$:
\begin{eqnarray}
&&\label{v1}\alpha_3(a,b\vdash c)=\alpha_3(a\vdash b,c)
=\alpha_3(b,a\vdash c)=\alpha_3(a, c\vdash b),\\
&&\alpha_2(a,b\vdash c)
=\alpha_2(b\vdash a,c),\\
&&\alpha_2(a, c\vdash b)=\alpha_2(b\vdash a,c)-\alpha_2(a\vdash b,c),\\
&&\alpha_1(a,b\vdash c)=\alpha_1(b\vdash a,c)=\alpha_1(c, b\vdash a),\\
&&\label{v2}\alpha_0(a,b\vdash c)-\alpha_0(b,a\vdash c)=\alpha_0(b\vdash a,c)+\alpha_0(b, c\vdash a).
\end{eqnarray}
Moreover,  such two 2-cocycles $\alpha_\lambda(\cdot,\cdot)$ and $\alpha_\lambda^{'}(\cdot,\cdot)$ are equivalent if and only if $\alpha_i(a,b)=\alpha_i^{'}(a,b)$ for $i=0$, $i=2$ and $i=3$, and there exists a linear map
$\varphi: V\rightarrow \mathbb{C}$ such that
\begin{eqnarray}
\alpha_1(a,b)=\alpha_1^{'}(a,b)+\varphi(a\vdash b+b\vdash a),
\end{eqnarray}
for all $a$, $b\in V$.
\end{coro}
\begin{proof}
It can be obtained from Theorem \ref{th4} with a similar proof as that in Corollary \ref{co6}.
\end{proof}
\begin{coro}
Let $(V,\vdash)$ be a Novikov algebra and $R=\mathbb{C}[\partial]V$ be the corresponding quadratic Lie conformal algebra. Set $\alpha_i(\cdot,\cdot)$ $(i=0, 1, 2, 3)$
be bilinear forms on $V$  satisfying (\ref{v1})-(\ref{v2}). Define bilinear forms $\varphi_i: \text{Coeff}(R)\times \text{Coeff}(R) \rightarrow \mathbb{C}$ as follows:
\begin{eqnarray}
&&\varphi_0(a\otimes t^m, b\otimes t^n)=\alpha_0(a,b)\delta_{m+n+1,0},\\
&&\varphi_1(a\otimes t^m, b\otimes t^n)=m\alpha_1(a,b)\delta_{m+n,0},\\
&&\varphi_2(a\otimes t^m, b\otimes t^n)=m(m-1)\alpha_2(a,b)\delta_{m+n-1,0},\\
&&\varphi_3(a\otimes t^m, b\otimes t^n)=m(m-1)(m-2)\alpha_3(a,b)\delta_{m+n-2,0},
\end{eqnarray}
for $a$, $b\in V$, $m$, $n\in \mathbb{Z}$. Then $\varphi_0$, $\varphi_1$, $\varphi_2$, $\varphi_3$ are 2-cocycles of Leibniz algebra $\text{Coeff}(R)$.
\end{coro}
\begin{proof}
The proof is similar to that in Corollary \ref{co7}.
\end{proof}

\begin{examp}
In this example, we give a characterization of Leibniz central extensions of $LC(C1)$, $LC(C2)$, $LC(C3(l))$ and
$LC(B)$ in Example \ref{exs1} by a one-dimensional center $\mathbb{C}\mathfrak{c}$ up to equivalence.

Note that Novikov algebras (C2), (C3(l)) and (B) satisfy the condition $V=V\vdash V$ in Corollary \ref{cp1}. Therefore, for determining the Leibniz central extensions of $LC(C2)$,  $LC(C3(l))$ and $LC(B)$, we only need to study $\alpha_i(\cdot,\cdot)$
for $i\in\{0,\cdots,3\}$ satisfying (\ref{v1})-(\ref{v2}) in the three cases.

First, we consider the case of $LC(C2)$. In this case, by a simple computation, we can get
that $\alpha_3(e_1,e_1)=\alpha_3(e_1,e_2)=\alpha_3(e_2,e_1)=0$, $\alpha_3(e_2,e_2)=\beta_1$,
$\alpha_2(\cdot,\cdot)=\alpha_0(\cdot,\cdot)=0$, $\alpha_1(e_1,e_2)=\alpha_1(e_2,e_1)=\beta_2$, and $\alpha_1(e_2,e_2)=\beta_3$, where $\beta_1$, $\beta_2$, $\beta_3\in \mathbb{C}$. By choosing the linear map $\varphi: V\rightarrow\mathbb{C}$ in Corollary \ref{cp1} defined by
$\varphi(e_1)=\beta_2$ and $\varphi(e_2)=\frac{\beta_3-2\beta_2}{2}$, we can make $\alpha_1(\cdot,\cdot)$ be zero up to equivalence. Therefore, by Corollary \ref{cp1}, all equivalence classes of central extensions of
$LC(C2)$ by a one-dimensional center $\mathbb{C}\mathfrak{c}$ are $\widehat{LC(C2)}(\beta_1)$  with the $\lambda$-product as follows:
\begin{eqnarray}
&&[{e_1}_\lambda e_1]=0,~~[{e_1}_\lambda e_2]=\lambda e_1,~~[{e_2}_\lambda e_1]=(\partial+\lambda)e_1,\\
&&[{e_2}_\lambda e_2]=(\partial+2\lambda)(e_1+e_2)+\beta_1\lambda^3\mathfrak{c},
\end{eqnarray}
for all $\beta_1\in \mathbb{C}$. Obviously, if $\beta_1\neq \beta_1^{'}$, then $\widehat{LC(C2)}(\beta_1)$ is not equivalent to $\widehat{LC(C2)}(\beta_1^{'})$. Therefore, $\text{dim}~H^2(LC(C2),\mathbb{C}\mathfrak{c})=1$.

With a similar discussion as above, we can get: \\
If $l\neq -1$ and $l\neq -2$, all equivalence classes of central extensions of
$LC(C3(l))$ by a one-dimensional center $\mathbb{C}\mathfrak{c}$ are $\widehat{LC(C3(l))}(\beta_1)$  with the $\lambda$-product as follows:
\begin{eqnarray}
&&[{e_1}_\lambda e_1]=0,~~[{e_1}_\lambda e_2]=(l\partial+(l+1)\lambda) e_1,~~[{e_2}_\lambda e_1]=(\partial+(l+1)\lambda)e_1,\\
&&[{e_2}_\lambda e_2]=(\partial+2\lambda)e_2+\beta_1\lambda^3\mathfrak{c},
\end{eqnarray}
for all $\beta_1\in \mathbb{C}$. Obviously, if $\beta_1\neq \beta_1^{'}$, then $\widehat{LC(C3(l))}(\beta_1)$ is not equivalent to $\widehat{LC(C3(l))}(\beta_1^{'})$. Therefore, if $l\neq -1$ and $l\neq -2$, then $\text{dim}~H^2(LC(C3(l)),\mathbb{C}\mathfrak{c})=1$;\\
All equivalence classes of central extensions of
$LC(C3(-1))$ by a one-dimensional center $\mathbb{C}\mathfrak{c}$ are $\widehat{LC(C3(-1))}(\beta_1,\beta_2, \beta_3)$  with the $\lambda$-product as follows:
\begin{eqnarray}
&&[{e_1}_\lambda e_1]=0,~~[{e_1}_\lambda e_2]=-\partial e_1+(-\beta_3+\beta_2\lambda)\mathfrak{c},\\
&&[{e_2}_\lambda e_1]=\partial e_1+(\beta_3+\beta_2\lambda)\mathfrak{c},\\
&&[{e_2}_\lambda e_2]=(\partial+2\lambda)e_2+\beta_1\lambda^3\mathfrak{c},
\end{eqnarray}
for all $\beta_1$, $\beta_2$, $\beta_3\in \mathbb{C}$. Obviously, if $(\beta_1,\beta_2, \beta_3)\neq (\beta_1^{'}, \beta_2^{'}, \beta_3^{'})$, then $\widehat{LC(C3(-1))}(\beta_1,\beta_2, \beta_3)$ is not equivalent to $\widehat{LC(C3(-1))}(\beta_1^{'}, \beta_2^{'}, \beta_3^{'})$. Therefore, $\text{dim}~H^2(LC(C3(-1)),\mathbb{C}\mathfrak{c})=3$;\\
All equivalence classes of central extensions of
$LC(C3(-2))$ by a one-dimensional center $\mathbb{C}\mathfrak{c}$ are $\widehat{LC(C3(-2))}(\beta_1,\beta_2)$  with the $\lambda$-product as follows:
\begin{eqnarray}
&&[{e_1}_\lambda e_1]=0,~~[{e_1}_\lambda e_2]=(-2\partial-\lambda) e_1+\beta_2\mathfrak{c},\\
&&[{e_2}_\lambda e_1]=(\partial-\lambda)e_1,\\
&&[{e_2}_\lambda e_2]=(\partial+2\lambda)e_2+\beta_1\lambda^3\mathfrak{c},
\end{eqnarray}
for all $\beta_1$, $\beta_2\in \mathbb{C}$. Obviously, if $(\beta_1,\beta_2)\neq (\beta_1^{'}, \beta_2^{'})$, then $\widehat{LC(C3(-2))}(\beta_1,\beta_2)$ is not equivalent to $\widehat{LC(C3(-2))}(\beta_1^{'}, \beta_2^{'})$. Therefore, $\text{dim}~H^2(LC(C3(-2)),\mathbb{C}\mathfrak{c})=2$. Note that if $\beta_2\neq 0$, then $\widehat{LC(C3(-2))}(\beta_1,\beta_2)$ is not a Lie conformal algebra.\\
All equivalence classes of central extensions of
$LC(B)$ by a one-dimensional center $\mathbb{C}\mathfrak{c}$ are $\widehat{LC(B)}(\beta_1,$ $\beta_2, \beta_3)$  with the $\lambda$-product as follows:
\begin{eqnarray}
&&[{e_1}_\lambda e_1]=\beta_3\lambda\mathfrak{c},~~[{e_1}_\lambda e_2]=\lambda e_1+\beta_2\lambda^2\mathfrak{c},\\
&&[{e_2}_\lambda e_1]=(\partial+\lambda)e_1-\beta_2\lambda^2\mathfrak{c},\\
&&[{e_2}_\lambda e_2]=(\partial+2\lambda)e_2+\beta_1\lambda^3\mathfrak{c},
\end{eqnarray}
for all $\beta_1$, $\beta_2$, $\beta_3\in \mathbb{C}$. Obviously, if $(\beta_1,\beta_2, \beta_3)\neq (\beta_1^{'}, \beta_2^{'}, \beta_3^{'})$, then $\widehat{LC(B)}(\beta_1,\beta_2, \beta_3)$ is not equivalent to $\widehat{LC(B)}(\beta_1^{'}, \beta_2^{'}, \beta_3^{'})$. Therefore, $\text{dim}~H^2(LC(B),\mathbb{C}\mathfrak{c})=3$.

Finally, we determine the central extensions of $LC(C1)$ by a one-dimensional center $\mathbb{C}\mathfrak{c}$ up to equivalence. By the conformal cocycle condition (\ref{cy4}) with $[\cdot,\cdot]$ trivial, it is easy to see that
$\text{dim}~H^2(LC(C1),\mathbb{C}\mathfrak{c})=\infty$.
\end{examp}

{\bf Acknowledgments}
{We wish to thank the referee for careful reading and useful comments. In particular, we would like to thank the referee
for suggesting the improvements of the present form of this paper.}


\small
\begin{thebibliography}{9999}\vskip0pt\small
\def\re{\bibitem}\parindent=2ex\parskip=-2pt\baselineskip=-2pt

\bibitem[1]{B} A.~Bloh, On a generalization of the concept of Lie algebra (Russian), Dokl.Akad.Nauk SSSR, {\bf 165} (1965), 471-473.
\bibitem[2]{BDK1} B.~Bakalov, A.~D'Andrea, V.~Kac, Theory of finite pseudoalgebras, Adv.Math. {\bf 162} (2001), 1-140.
\bibitem[3]{BDK} A.Barakat, A.De sole, V.Kac, Poisson vertex algebras in the theory of Hamiltonian equations, Japan. J. Math. {\bf 4}(2009), 141-252.
\bibitem[4]{BKV} B.~Bakalov, V.~Kac, A.~Voronov,  Cohomology of conformal algebras, Comm. Math. Phys.  {\bf 200} (1999),561-598.
\bibitem[5]{BM}C.~M. Bai, D.~J. Meng, The classification of Novikov algebras in low dimensions, J. Phys. A: Math. Gen. {\bf 34} (2001): 1581-1594.
\bibitem[6]{BN}  A.~A. Balinskii, S.~P. Novikov,
Poisson brackets of hydrodynamical type, Frobenius algebras and Lie
algebras,  Dokladu AN SSSR, {\bf 283} (1985),1036-1039.

\bibitem[7]{Ca}  A.~Cayley,  On the Theory of Analytic Forms Called Trees,
Collected Mathematical Papers of Arthur Cayley, Cambridge Univ.
Press, Cambridge, {\bf 3} (1890), 242-246.
\bibitem[8]{Ch1}  F.~Chapoton, Un endofoncteur de la cat$\acute{e}$gorie des op$\acute{e}$rades, In: Loday J.-L., Frabetti A., Chapoton F., Goichot F. (Eds), \emph{Dialgeras and related operads}, Spring-Verl., Berlin, 2001, pp. 105-110. (Lectures Notes in Math., vol.1763).
\bibitem[9]{Ch} B.~Y. Chu, Symplectic homogeneous spaces, Trans. Amer. Math. Soc. {\bf 197} (1974),145-159.
\bibitem[10]{CK} A.~Connes, D.~Kreimer, Hopf algebras, renormalization and noncommutative geometry, Comm. Math. Phys. {\bf 199} (1998),203-242.
\bibitem[11]{CK1} S.~Cheng, V.~Kac,  Conformal modules, Asian J.Math. {\bf 1}(1997), 181-193.
\bibitem[12]{CK2} S.~Cheng, V.~Kac, M.~Wakimoto,  Extensions of conformal modules, in:Topological Field Theory, Primitive Forms and Related Topics(Kyoto), in: Progress in Math., Vol.160, Birkh$\ddot{a}$user, Boston, 1998, pp. 33-57; q-alg/9709019.
\bibitem[13]{DK} A.~D'Andrea, V.~Kac, Structure theory of
finite conformal algebras, Selecta Math. {\bf 4} (1998), 377-418.
\bibitem[14]{GD} I.~M. Gel'fand, I.~Ya. Dorfman, Hamiltonian operators and algebraic structures related to
them, Funkts.Anal.Prilozhen {\bf 13}(1979), 13-30.
\bibitem[15]{gk}V.~Yu. Gubarev, P.~S. Kolesnikov, Operads of decorated trees and their duals, Comment. Math. Univ. Carolin. {\bf 55} (2014), 421-445.
\bibitem[16]{gk1}V.~Yu. Gubarev, P.~S. Kolesnikov, Embedding of dendriform algebras into Rota-Baxter algebras, Centr. Eur. J. Math. {\bf 11} (2013), 226-245.
\bibitem[17]{GP}L.~Guo, J.~Pei, Averaging algebras, Schr\"{o}der numbers, and rooted trees, J.Algebraic Combin. {\bf 42}(2015), 73-109.
\bibitem[18]{H1}Y.~Y. Hong, Central extensions and conformal derivations of a
 class of Lie conformal algebras, arXiv:1502.02770v4.
\bibitem[19]{K1} V.~Kac, Vertex Algebras for Beginners, 2nd Edition, Amer. Math. Soc., Providence, RI,1998.
\bibitem[20]{K2} V.~Kac, Lie superalgebras, Adv.Math. {\bf 26} (1977), 8-96.
\bibitem[21]{L}J.-L.Loday, Cyclic Homology, Grundlehren der Mathematischen Wissenschaften, Springer-Varlag, Berlin, Vol. 301, 1991.
\bibitem[22]{Os}  J.~M.~Osborn, Novikov algebras, Nova J. Algebra \& Geom. {\bf 1} (1992), 1-14.
\bibitem[23]{Xu1}  X.~P. Xu, Quadratic conformal superalgebras, J.Algebra {\bf 231}(2000), 1-38.
\bibitem[24]{Wu} Z.~X. Wu, Leibniz $H$-pseudoalgebras, J.Algebra, {\bf 437}(2015), 1-33.
\bibitem[25]{Z}J. Zhang, On the cohomology of Leibniz conformal algebras, J.Math.Phys. {\bf 56} (2015), 041703.
\end{thebibliography}
\end{document}